\documentclass{commat}

\title{%
    Legendre curves on 3-dimensional  $C_{12}$-Manifolds
    }

\author{%
    Gherici Beldjilali, Benaoumeur Bayour and Habib Bouzir
    }

\affiliation{
    \address{Gherici Beldjilali --
    Laboratory of Quantum Physics and Mathematical Modeling (LPQ3M)\\
    University of Mascara,   Algeria.
        }
    \email{%
    gherici.beldjilali@univ-mascara.dz
    }
    \address{Benaoumeur Bayour --
    Departement of  Mathematics, University of Mascara, Algeria.
        }
    \email{%
    b.bayour@univ-mascara.dz
    }
    \address{Habib Bouzir --
    Laboratory of Quantum Physics and Mathematical Modeling (LPQ3M)\\
    University of Mascara,   Algeria.
        }
    \email{%
    habib.bouzir@univ-mascara.dz
    }
    }

\abstract{%
   Legendre curves play a very important and special role in geometry and topology of almost contact manifolds.  There are certain results known for Legendre curves in 3-dimensional normal almost contact manifolds.  The aim of this paper is to study Legendre curves of three-dimensional $C_{12}$-manifolds which are non-normal almost contact manifolds and classifying all biharmonic Legendre curves in these manifolds.
    }

\keywords{%
    Almost contact metric structure, $C_{12}$-manifolds, Legendre curves.
    }

\msc{%
    Primary 53D15, Secondary 53B25.
    }

\VOLUME{31}
\YEAR{2023}
\NUMBER{1}
\firstpage{293}
\DOI{https://doi.org/10.46298/cm.10390}

\begin{paper} 

\section{Introduction}
In \cite{CG}, D. Chinea and C. Gonzalez have defined 12 classes of  almost contact metric manifolds. In dimension 3, these manifolds are reduced to five classes: $\vert C \vert$ class of cosymplectic manifolds, $C_5$ class of $\beta$-Kenmotsu manifolds, $C_6$ class of $\alpha$-Sasakian manifolds, $C_9$-manifolds and $C_{12}$-manifolds.

Only the last two classes can never be normal. For this reason, all work concerning curves on almost contact metric manifolds focuses on the first three classes. for example, Legendre curves on contact manifolds have been studied by C. Baikoussis and
D.E. Blair in the paper \cite{BB}. M. Belkhelfa et al. \cite{BH} have investigated Legendre curves in Riemannian and
Lorentzian manifolds. Also, Legendre curves have been studied on Kenmotsu manifolds in \cite{TMZC}, on $\alpha$-Sasaki minifolds in \cite{OT}, on quasi-Sasakian manifolds in \cite{WEL}, on Trans-Sasaki manifolds in \cite{SHM} and others.  

$C_{12}$-manifolds arose in a natural way from the classification of almost contact metric structures by D. Chinea and C. Gonzalez \cite{CG} in 1990. Recently, these manifolds was studied and most of its properties were deduced in \cite{CF}, especially \cite{BBB}. As $C_{12}$-manifolds differs from the above almost contact metric  manifolds, at least  in the normality property, in the present paper, we will show that the results related to Legendre curves in 3-dimentional $C_{12}$-manifolds are fundamentally different from the results in other  almost contact metric  manifolds.

The present paper is organized as follows: After the introduction, we recall   some required preliminaries on 
  almost contact geometry and Frenet curves in general Riemannian geometry in Section 3. In Section 4, we give a
very brief review of 3-dimentional $C_{12}$-manifolds with concrete example. The next section is focused on the study of Legendre curves on  3-dimentional $C_{12}$-manifolds.  In the last section, we  demonstrate a nice property of  biharmonic Legendre curves in 3-dimentional $C_{12}$-manifolds.

\section{Almost contact manifold}
An odd-dimensional Riemannian manifold $(M^{2n+1},g)$ is said to be an almost contact metric manifold if there exists on $M$ a $(1,1)$-tensor field $\varphi$, a vector field $\xi$ (called the structure vector field) and a $1$-form $\eta$ such that 
\begin{equation}\label{CondPresqCont}
 \left\{
          \begin{array}{lll}
\eta(\xi)=1,\\
\varphi^{2}(X) = -X+\eta(X)\xi ,\\
g(\varphi X,\varphi Y) = g(X,Y)-\eta(X)\eta(Y),
          \end{array}
  \right.
\end{equation}
 for any vector fields $X$, $Y$ on $M$. 
 
 In particular, in an almost contact metric manifold we also have
 $$ \varphi\xi=0 \qquad  and \qquad \eta \circ \varphi=0.$$

The fundamental 2-form $\phi$ is defined by 
$$\phi (X,Y) = g(X, \varphi Y).$$
It is known that the almost contact
structure $(\varphi, \xi, \eta)$  is said to be normal if and only if
\begin{equation}\label{N1phi}
N^{(1)}(X,Y)=N_{\varphi}(X,Y)+2d\eta(X,Y)\xi=0,
\end{equation}
 for any $X$, $Y$ on $M$, where $N_{\varphi}$ denotes the Nijenhuis torsion of $\varphi$, given by
\begin{equation}\label{Nphi}
N_{\varphi}(X,Y)=\varphi^{2}[X,Y]+[\varphi X,\varphi Y]-\varphi[\varphi X,Y]-\varphi[X,\varphi Y].
\end{equation}

 For more background on almost contact metric manifolds, we recommend the references \cite{BLA}, \cite{BGM}, \cite{YK}.

\section{Legendre curves}

Let $(M, g)$ be a 3-dimensional Riemannian manifold with Levi-Civita connection $\nabla$ and $\gamma: I \rightarrow M$ parameterized by the arc length. $\gamma$ is  said  to  be  a Frenet curve if  there  exists  an  orthonormal frame $\{ E_1= \dot{\gamma}, E_2, E_3 \}$ along $\gamma$ such that
\begin{equation}\label{RFrenet}
 \left\{
          \begin{array}{lll}
 \nabla_{\dot{\gamma}} E_1 = \kappa E_2,\\
\nabla_{\dot{\gamma}} E_2 = -\kappa E_1 + \tau E_3,\\
\nabla_{\dot{\gamma}} E_3 = -\tau E_2.
          \end{array}
  \right.
\end{equation}
The curvature $\kappa$ is defined by the formula 
\begin{equation}\label{Kappa}
 \kappa = \vert  \nabla_{\dot{\gamma}} \dot{\gamma}\vert.
 \end{equation}
The second unit vector field $E_2$ is thus obtained by
\begin{equation}\label{2unitVector}
  \nabla_{\dot{\gamma}} \dot{\gamma}= \kappa E_2.
   \end{equation}
Next, the torsion $\tau$ and the third unit vector field $E_3$ are defined by the formulas
\begin{equation}\label{3unitVector}
 \tau = \vert  \nabla_{\dot{\gamma}} E_2 + \kappa E_1 \vert \qquad and \qquad
\nabla_{\dot{\gamma}} E_2 + \kappa E_1 = \tau E_3.
 \end{equation}
 
A Frenet curve $\gamma: I \rightarrow M$ in an almost contact metric manifold is said to be a Legendre curve \cite{BB}, if it is an integral curve of the contact distribution $\mathcal{D} = ker\; \eta$. Formally, it is also said that a Frenet curve $\gamma$ in an almost contact metric manifold is a Legendre curve if and only if $\eta(\dot{\gamma})=0$ and $g(\dot{\gamma} , \dot{\gamma}) = 1$. 

\section{Three dimensional  $C_{12}$-manifold} 

In the classification of D. Chinea and C. Gonzalez \cite{CG} of almost contact metric manifolds there is a class $C_{12}$-manifolds  which can be integrable but never normal. In this classification, $C_{12}$-manifolds are defined by
\begin{equation}\label{CGC12}
(\nabla_X \phi)(Y,Z)=\eta(X)\eta(Z) (\nabla_{\xi} \eta) \varphi Y - \eta(X)\eta(Y) (\nabla_{\xi} \eta) \varphi Z.
\end{equation}

In \cite{BBB} and \cite{CF}, The $(2n+1)$-dimensional $C_{12}$-manifolds is characterized by:

\begin{equation}\label{NablaVarphi}
 (\nabla_X \varphi) Y =  \eta(X) \big( \omega(\varphi Y) \xi + \eta(Y) \varphi \psi \big),
\end{equation}
for any $X$ and $Y$ vector fields on $M$, where  $ \omega = -\big( \nabla_{\xi}\xi \big)^{\flat} = -\nabla_{\xi}\eta$ and $\psi$ is the vector field given by 
$$\omega(X)=g(X, \psi) = -g(X, \nabla_{\xi}\xi), $$ for all $X$ vector field on $M$.

Moreover,  in  \cite{BBB} the $(2n+1)$-dimensional $C_{12}$-manifolds is also  characterized by 
\begin{equation}\label{cond11}
 {\rm d} \eta =   \omega \wedge \eta  \qquad  {\rm d} \phi = 0 \qquad  and \qquad N_{\varphi} =0.
\end{equation}

Here, we emphasize that the almost $C_{12}$-manifolds is defined by the following:
\begin{definition}
Let $(M^{2n+1}, \varphi, \xi, \eta,g)$ be an almost contact manifold. $M$ is called almost $C_{12}$-manifold if there exists a closed one-form $\omega$  which satisfies
\begin{equation*}
 {\rm d} \eta =   \omega \wedge \eta  \qquad and \qquad {\rm d} \phi = 0.
\end{equation*}
In addition, if $N_{\varphi} =0$ we say that $M$ is a $C_{12}$-manifold.
\end{definition}

In \cite{BBB}, the authors studied the 3-dimensional unit $C_{12}$-manifold i.e. the case where $\psi$ is a unit vector field. We will deal here with the general case, i.e. $\psi$ is not necessarily unitary. For that, taking $V = {\rm e}^{-\rho}\psi$ where $ {\rm e}^{\rho} = \vert \psi \vert$, we get immediately that  $\{\xi, V, \varphi V \}$ is  an orthonormal frame. We refer to this basis as  Fundamental basis.

Using this frame, one can get the following:
\begin{proposition}\label{prop3}
For any $C_{12}$-manifold, for all vector field $X$ on $M$ we have
\begin{itemize}
  \item[1)] $\nabla_{X} \xi =  - {\rm e}^{\rho}\eta(X) V$
  \item[2)] $\nabla_{\xi }V= {\rm e}^{\rho} \xi$
  \item[3)] $ \nabla_V  V = \varphi V(\rho) \varphi V$
  \item[4)] $ \nabla_\xi \varphi V = 0$
  \item[5)] $ \nabla_V  \varphi V =-\varphi V(\rho) V. $ 
 \end{itemize}
\end{proposition}
\begin{proof}
For the first, using (\ref{NablaVarphi}) for $Y=\xi$ we get
\begin{eqnarray*}
(\nabla_X \varphi)\xi &=& \eta(X) \varphi \psi \\
&=& {\rm e}^{\rho}\eta(X) \varphi V,
\end{eqnarray*}
knowing that $ (\nabla_X \varphi )Y =  \nabla_X \varphi Y - \varphi \nabla_X Y $ and applying $\varphi$ we obtain
\begin{eqnarray*}
\nabla_X \xi &=& {\rm e}^{\rho}\eta(X) \varphi^2 V\\
&=& -{\rm e}^{\rho}\eta(X)  V.
\end{eqnarray*}
 For the second, we have 
\begin{eqnarray*}
2 {\rm d}\omega(\xi , X) = 0 \Leftrightarrow  g(\nabla_{\xi} \psi , X) &=&g( \nabla_X \psi, \xi) \\
&=& - g(\psi, \nabla_X \xi )\\
&=&  {\rm e}^{2\rho}\eta(X),
\end{eqnarray*}
 which gives $ \nabla_{\xi} \psi={\rm e}^{2\rho} \xi$ and  then
 \begin{eqnarray*}
\nabla_{\xi }V &=&  \nabla_{\xi } ({\rm e}^{-\rho}\psi) \\
&=&- \xi(\rho)V + {\rm e}^{\rho} \xi .
\end{eqnarray*} 
 On the other hand, we have
  \begin{eqnarray*}
\xi(\rho) &=&  \frac{1}{2} {\rm e}^{-2\rho}\xi({\rm e}^{2\rho}) \\
&=& \frac{1}{2} {\rm e}^{-2\rho}\xi \big( g(\psi , \psi) \big)\\
&=& {\rm e}^{-2\rho} g(\nabla_{\xi}\psi , \psi)=0,
\end{eqnarray*} 
 then, 
 $$ \nabla_{\xi }V ={\rm e}^{\rho} \xi.$$
For $\nabla_{V} V$, we have 
\begin{eqnarray*}
2 {\rm d}\omega(\psi , X) = 0 \Leftrightarrow  g(\nabla_{\psi} \psi , X) &=&g( \nabla_X \psi, \psi) \\
&=& \frac{1}{2}Xg(\psi, \psi)\\
&=& {\rm e}^{2\rho} g( {\rm grad}  \rho, X),
\end{eqnarray*}
i.e. $\nabla_{\psi} \psi ={\rm e}^{2\rho} {\rm grad}  \rho$ which gives $ \nabla_V V = {\rm grad}  \rho -V(\rho) V.$\\
 Also, we have
\begin{eqnarray*}
{\rm grad}  \rho&=& \xi(\rho) \xi + V(\rho) V + \varphi V(\rho) \varphi V\\
&=& V(\rho) V + \varphi V(\rho) \varphi V,
\end{eqnarray*}
then,
$$ \nabla_V V = \varphi V(\rho) \varphi V.$$
For the rest, just use the formula $\nabla_X \varphi Y =( \nabla_X \varphi )Y + \varphi \nabla_X Y$ noting that
\[ (\nabla_{V}\varphi )X = (\nabla_{\varphi V}\varphi )X=0. \qedhere\]
\end{proof}
It remains to count $\nabla_{\varphi V}V$ and $\nabla_{\varphi V} \varphi V$.  For that, we have the following lemma
\begin{lemma}\label{Lem4}
For any $3$-dimensional  $C_{12}$-manifold, we have
\begin{itemize}
  \item[1)] $\nabla_{\varphi V} V =( -{\rm e}^{\rho} + {\rm div} V ) \varphi V,$
  \item[2)] $\nabla_{\varphi V } \varphi V =( {\rm e}^{\rho} - {\rm div} V )V.$
 \end{itemize}
\end{lemma}
\begin{proof}
Since $\{\xi, V, \varphi V \}$ is  an orthonormal frame then,
$$ \nabla_{\varphi V } V= a\; \xi + b\;  V +  c\; \varphi V,$$
Using Proposition \ref{prop3}, we have
 $$a=g( \nabla_{\varphi V} V , \xi) = -g( V , \nabla_{\varphi V } \xi ) = 0$$
  and  
$b = g(  \nabla_{\varphi V} V , V) =0 $. for get $c$ we have
\begin{eqnarray*}
{\rm div} V &=&g( \nabla_{\xi } V, \xi) + g( \nabla_{\varphi V  } \psi , \varphi V ) \\
&=&{\rm e}^{\rho} +  g( \nabla_{\varphi \psi  } \psi , \varphi \psi )
\Leftrightarrow  g( \nabla_{\varphi V  } V , \varphi V) = -{\rm e}^{\rho} + {\rm div} V,
\end{eqnarray*}
then,
$$\nabla_{\varphi V} V =( -{\rm e}^{\rho} + {\rm div} V ) \varphi V.$$
Applying $\varphi$ with (\ref{NablaVarphi}), we obtain
\[\nabla_{\varphi V } \varphi V =( {\rm e}^{\rho} - {\rm div} V )V.\qedhere\]
\end{proof}

According to the Proposition \ref{prop3} and Lemma \ref{Lem4},  the 3-dimensional $C_{12}$-manifold is completely controllable. That is:
\begin{corollary}\label{Cor2}
For any $C_{12}$-manifold, we have
\begin{equation*}
          \begin{array}{lll}
  \nabla_{\xi} \xi =  - {\rm e}^{\rho}V,
 &  \nabla_{\xi }V= {\rm e}^{\rho} \xi,
 & \nabla_\xi \varphi V =0,\\
 \nabla_{V} \xi =0,
 &\nabla_V  V = \varphi V(\rho) \varphi V,
 &  \nabla_{V} \varphi V =-\varphi V(\rho) V,\\
\nabla_{\varphi V} \xi =0,
& \nabla_{\varphi V} V =( -{\rm e}^{\rho} + {\rm div} V ) \varphi V,
& \nabla_{\varphi V } \varphi V =( {\rm e}^{\rho} - {\rm div} V )V.
          \end{array}
\end{equation*}
\end{corollary}

\begin{example}\label{Ex1}
We denote the Cartesian coordinates in a $3$-dimensional Euclidean space $M=\mathbb{R}^3$ by $(x,y,z)$ and define a symmetric tensor field $g$ by
\begin{eqnarray*}
g= {\rm e}^{2y}\left(
       \begin{array}{ccc}
        1+\alpha^2 & 0 & -1\\
        0 & \alpha^2 & 0\\
        -1 & 0 & 1
       \end{array}
\right),        
\end{eqnarray*}
where $\alpha=\alpha(x,y) \geq 0$ every where is a function on $\mathbb{R}^3$.
Further, we define an almost contact metric  $(\varphi ,\xi ,\eta)$ on $\mathbb{R}^3$ by
$$
\varphi= \left(
       \begin{array}{ccc}
        0 & -1 & 0\\
        1 &  0 & 0\\
        0 & -1& 0
       \end{array}
\right),\quad
   \xi= {\rm e}^{-y}\left(
       \begin{array}{c}
        0\\
        0\\
        1
       \end{array}
\right),\quad
\begin{array}{lll}
  \eta = {\rm e}^{y}(-1,0,1).
\end{array}
$$
The fundamental 1-form $\eta$ and the 2-form  $\phi$ have the forms,
$$\eta= {\rm e}^{y} (d z- d x) \qquad and \qquad \phi= -2 \alpha^2 {\rm e}^{2y}d x \wedge d y,$$
and hence
$$
{\rm d}\eta = {\rm e}^{y} \Big(   d x\wedge d y+  d y\wedge  d z \Big) =  dy \wedge \eta,$$
$$ {\rm d}\phi = 0.$$

By a direct computation the non trivial components of 
$N_{kj}^{(1)\;i}$ are given by
$$N_{12}^{(1)\;3}=1 , \quad N_{23}^{(1)\;3}=1.$$
i.e. $ N^{(1)}\neq 0$. But,  $ \forall i,j,k \in \{1,2,3\}$ 
\begin{eqnarray*}
(N_{\varphi})_{kj}^{i}=0,
\end{eqnarray*}
implying that  $(\varphi, \xi, \eta)$ becomes integrable non normal. We have  $\omega =  dy$  i.e.  ${ \rm d} \omega =0$
and knowing that $\omega$ is the $g$-dual of $\psi$ i.e. $\omega (X)=g(X, \psi)$, we have immediately that
\begin{equation}\label{Psi}
\psi=  \frac{{\rm e}^{-2y}}{\alpha^2}  \frac{\partial}{\partial y}.
\end{equation}
Thus, $(\varphi, \xi, \psi, \eta, \omega, g)$ is a 1-parameter family of $C_{12}$-structure on $\mathbb{R}^3$.\\
Notice that 
$$\vert \psi \vert^2 =\omega(\psi)= g(\psi , \psi) = \frac{{\rm e}^{-2y}}{\alpha^2} $$
implies $V = \frac{{\rm e}^{-y}}{\alpha} \frac{\partial}{\partial y}$ is a unit vector field, then 
$$\left\lbrace \xi = {\rm e}^{-y} \frac{\partial}{\partial z}, \quad V = \frac{{\rm e}^{-y}}{\alpha} \frac{\partial}{\partial y},
\quad \varphi V =\frac{{\rm e}^{-y}}{\alpha}  \big( \frac{\partial}{\partial x} +  \frac{\partial}{\partial z}\big)\right\rbrace $$
   form an orthonormal basis.
To verify result in formula (\ref{NablaVarphi}), the  components of the Levi-Civita connection corresponding to $g$ are given by:
\begin{equation*}
          \begin{array}{lll}
 \nabla_{\xi} \xi =-\frac{{\rm e}^{-y}}{\alpha} V,
 &  \nabla_{\xi} V =\frac{{\rm e}^{-y}}{\alpha}  \xi,
 & \nabla_{\xi} \varphi V =0,\\
 \nabla_{V} \xi =0,
 & \nabla_{V} V =- \frac{{\rm e}^{-y}}{\alpha^2}  \alpha_1  \varphi V,
 &  \nabla_{V} \varphi V =- \varphi \nabla_{V} V ,\\
\nabla_{\varphi V} \xi =0,
& \nabla_{\varphi V} V =\frac{{\rm e}^{-y}}{\alpha^2} ( \alpha + \alpha_2 ) \varphi V,
&  \nabla_{\varphi V} \varphi V = \varphi \nabla_{\varphi V} V ,
          \end{array}
\end{equation*}
where $ \alpha_2 = \frac{\partial \alpha}{\partial y}$.  Then, one can easily check that for all $ i, j \in \{1,2,3\}$
 \begin{eqnarray*}
 (\nabla_{e_i} \varphi) e_j &=& \nabla_{e_i} \varphi e_j  - \varphi \nabla_{e_i}  e_j \\
   &=& \eta(e_i) \big( \omega(\varphi e_j) \xi +\eta(e_j) \varphi \psi \big).
\end{eqnarray*}
\end{example}

Through the rest of this paper $(M, \varphi, \xi, \psi, \eta, \omega, g)$ always denotes a 3-dimensional $C_{12}$-manifold and $\{\xi, V, \varphi V\}$ it's fundamental frame.

\section{Legendre curves on 3-dimensional  $C_{12}$-manifold}

Let $\gamma(s)$ be a Legendre curve in a  3-dimensional $C_{12}$-manifold $M$. 
Let us compute the equations of motion for the associated Frenet frame is $\{ \dot{\gamma}, \varphi  \dot{\gamma}, \xi\}$.

 Differentiating $\eta(\dot{\gamma})= 0$ along $\gamma$ we get
$$   g(  \nabla_{\dot{\gamma}}\dot{\gamma},\xi)=-g(\nabla_{\dot{\gamma}} \xi , \dot{\gamma}),$$
with the help of Proposition \ref{prop3}, we get $g(\nabla_{\dot{\gamma}}\dot{\gamma} , \xi ) =0$,
then,
\begin{equation}\label{EqE1}
 \nabla_{\dot{\gamma}} \dot{\gamma} =   \kappa \;\varphi \dot{\gamma},
\end{equation}
where $\kappa = \vert g(  \nabla_{\dot{\gamma}} \dot{\gamma}  , \varphi \dot{\gamma} )\varphi \dot{\gamma} \vert$.\\
For $\nabla_{\dot{\gamma}} \varphi \dot{\gamma}$, using (\ref{NablaVarphi}) we obtain
\begin{eqnarray}\label{EqE2}
 \nabla_{\dot{\gamma}} \varphi \dot{\gamma} & =&   ( \nabla_{\dot{\gamma}} \varphi) \dot{\gamma} +   \varphi  \nabla_{\dot{\gamma}}\dot{\gamma}\notag \\
 & =&   \kappa \varphi^2 \dot{\gamma} \notag \\
  &=& - \kappa \dot{\gamma}.
\end{eqnarray}
Using Proposition \ref{prop3}, we get
\begin{equation}\label{EqE3}
 \nabla_{\dot{\gamma}} \xi =0.
\end{equation}
We compare the equations (\ref{EqE1})-(\ref{EqE3}) and (\ref{RFrenet}), we conclude the following:
 \begin{proposition}
In a 3-dimensional $C_{12}$-manifold, a Legendre curve is a plane Frenet curve.
\end{proposition}
Conversely, suppose that $\gamma$ is a plane Frenet curve in a 3-dimensional $C_{12}$-manifold. That is $\tau =0$
and  the equations of motion (\ref{RFrenet}) becomes
\begin{equation}\label{RFrenetTauZero}
 \left\{
          \begin{array}{lll}
 \nabla_{\dot{\gamma}} E_1 = \kappa E_2,\\
\nabla_{\dot{\gamma}} E_2 = -\kappa E_1,\\
\nabla_{\dot{\gamma}} E_3 =0.
          \end{array}
  \right.
\end{equation}
  Putting $\eta( \dot{\gamma})=\sigma$ then, since $\dot{\gamma}$ is different from $\xi$ and both of them are
unit vectors we observe that $0 \leq \vert g( \xi, \dot{\gamma})=\sigma \vert < 1$ on $M$. We
note that $\dot{\gamma}$ is not collinear with $\xi$. So, it can be verified that the vector fields
\begin{equation}\label{RepAss}   
E_1 = \dot{\gamma},\qquad E_2 = \frac{\varphi \dot{\gamma}}{\sqrt{1-\sigma^2}},
\qquad E_3 = \frac{\xi - \sigma \dot{\gamma}}{\sqrt{1-\sigma^2}},
\end{equation}
 form an orthonormal frame along $\gamma$ and consequently, we can write
 \begin{equation}\label{EqPsi} 
  \psi = \omega(\dot{\gamma}) E_1 + \frac{\omega\varphi \dot{\gamma})}{\sqrt{1-\sigma^2}}E_2 - \frac{\sigma \omega(\dot{\gamma})}{\sqrt{1-\sigma^2}}E_3.
  \end{equation}
 From (\ref{RepAss}) and (\ref{EqPsi}), one can get
  \begin{eqnarray*}
 \nabla_{\dot{\gamma}} E_3 &=& \nabla_{\dot{\gamma}} \Big(  \frac{\xi - \sigma \dot{\gamma}}{\sqrt{1-\sigma^2}}\Big) \\
  &=&  \frac{ -1}{\sqrt{1-\sigma^2}}\Big( \dot{\sigma} + \sigma \omega(\dot{\gamma}) \Big) E_1
            - \frac{\sigma }{1-\sigma^2}\Big( \omega( \varphi \dot{\gamma} ) + \kappa \sqrt{1-\sigma^2} \Big) E_2\\
          &&  +  \frac{\sigma }{1-\sigma^2}\Big(\dot{\sigma} + \sigma  \omega(\dot{\gamma} )  \Big) E_3.
  \end{eqnarray*}
  Since $\nabla_{\dot{\gamma}} E_3 =0$ then, we get
  \begin{equation}\label{Sys1}
 \left\{
          \begin{array}{ll}
 \dot{\sigma} + \sigma \omega(\dot{\gamma}) = 0\\
\sigma \big(\omega( \varphi \dot{\gamma} ) +  \kappa \sqrt{1-\sigma^2})=0.
          \end{array}
  \right.
\end{equation}
  This establishes the following theorem:
\begin{theorem}
 Let $M$ be a 3-dimensional $C_{12}$-manifold and  $ \gamma: I \rightarrow M$  a plane Frenet curve in $M$ (i.e. $\tau =0$), set  $\sigma = \eta(\dot{\gamma})$. If at a certain point $t_0 \in I$, $\sigma(t_0) = 0$, then $\gamma$ is a Legendre curve.
\end{theorem}

\section{Biharmonic curves on 3-dimensional  $C_{12}$-manifold} 

The purpose of this last section is to study the curves $ \gamma: I \rightarrow M$ which are biharmonic, i.e. which satisfy
$ \Delta^2 \gamma =0$ where $\Delta$  is the Laplacian of $I$ i.e. $\Delta = -\frac{d^2}{ds^2}$ with $s$ is the arclength parameter.

In \cite{CHE}, Chen defined a biharmonic submanifold $N \subset \mathbb{E}^n$  of the Euclidean
space as its mean curvature vector field $H$ satisfies $\Delta H = 0$.
 We note $ H = \nabla_{\dot{\gamma}}\dot{\gamma}$, then 
$$\Delta H = - \nabla_{\dot{\gamma}}\nabla_{\dot{\gamma}}H = 0.$$
Let $ \gamma$ be a Legendre curve 3-dimensional $C_{12}$-manifold parametrized by arc length and  its associated
Frenet frame is $\{\dot{\gamma}, \varphi \dot{\gamma} , \xi\}$. So, we have
$$  H = \nabla_{\dot{\gamma}}\dot{\gamma} = \kappa \varphi \dot{\gamma}.$$
With the help of (\ref{NablaVarphi}), we get
$$\nabla_{\dot{\gamma}} H =  \dot{\kappa } \varphi \dot{\gamma} - \kappa^2 \dot{\gamma},$$
then,
\begin{eqnarray*}
\Delta H &=& - \nabla_{\dot{\gamma}}\big(\dot{\kappa } \varphi \dot{\gamma} - \kappa^2 \dot{\gamma} \big)\\
&=& 3 \dot{\kappa }\kappa  \dot{\gamma} - (\ddot{\kappa } -\kappa^3) \varphi \dot{\gamma},
\end{eqnarray*}
therefore, $\Delta H = 0$ if and only if
\begin{equation}\label{Sys2}
 \left\{
          \begin{array}{ll}
 \dot{\kappa }\kappa = 0\\
\ddot{\kappa } -\kappa^3=0.
          \end{array}
  \right.
\end{equation}
Hence, $\kappa =0$ and we have the following result:
\begin{theorem}
 All biharmonic Legendre curves $\gamma$ in a 3-dimensional $C_{12}$-manifold are straight lines, i.e.
totally geodesic.
\end{theorem}


\EditInfo{November 12, 2021}{December 31, 2021}{Haizhong Li}

\end{paper}